\documentclass[11pt]{article}
\usepackage{amscd}
\usepackage{amsfonts}
\usepackage{amsmath}
\usepackage{amssymb}
\usepackage{amsthm}
\usepackage{bbm}
\usepackage{CJK}
\usepackage{fancyhdr}
\usepackage{graphicx}
\usepackage{hyperref}
\usepackage{indentfirst}
\usepackage{latexsym}
\usepackage{mathrsfs}
\usepackage{xypic}

\newtheorem{theorem}{Theorem}[section]
\newtheorem{lemma}[theorem]{Lemma}
\newtheorem{definition}[theorem]{Definition}

\newtheorem{corollary}[theorem]{Corollary}

\usepackage[top=1in,bottom=1in,left=1.25in,right=1.25in]{geometry}
\def\<{\langle}
\def\>{\rangle}
\def\a{\alpha}
\def\b{\beta}

\def\c{\cdot}

\def\g{\gamma}

\date{}
\begin{document}
\renewcommand{\baselinestretch}{1.2}
\renewcommand{\arraystretch}{1.0}
\title{\bf $n$-derivations of Lie color algebras}
\author{{\bf    Yizheng Li,  Shuangjian Guo\footnote
        { Corresponding author: shuangjianguo@126.com}}\\
{\small School of Mathematics and Statistics, Guizhou University of Finance and Economics} \\
{\small  Guiyang  550025, P. R. of China}}
 \maketitle

\begin{center}
\begin{minipage}{13.cm}

{\bf \begin{center} ABSTRACT \end{center}}
The aim of this article is to discuss the $n$-derivation algebras of Lie color algebras. It is proved that, if the
base ring contains $\frac{1}{n-1}$, $L$ is a perfect Lie color algebra with zero center, then every triple
derivation of $L$ is a derivation, and every $n$-derivation of the derivation algebra
$nDer(L))$ is an inner derivation.
\smallskip

\smallskip

 {\bf MR(2010) Subject Classification} 17B05, 17B40, 17B60
 \end{minipage}
 \end{center}
 \normalsize\vskip1cm

\section{INTRODUCTION}
\def\theequation{0. \arabic{equation}}
\setcounter{equation} {0}

The concept of derivations appeared in different mathematical fields with many
different forms. In algebra systems, derivations are linear maps that satisfy the Leibniz
relation. There are several kinds of derivations in the theory of  Lie algebras,
such as generalized derivations, Lie $n$-derivations,  $N$-derivations,  double derivations of
Lie algebras (\cite{Chen13}-\cite{Wang13}). In \cite{Zhou13}, Zhou studied triple derivations of Lie algebras.
It is proved that every triple derivation of a perfect Lie algebra with zero
center is a derivation. Moreover, every derivation of the derivation algebra is an
inner derivation.   Double derivations of Lie superalgebras were introduced in \cite{Sun18}, these
derivations are similar to the triple derivations of Lie algebras to some extent in \cite{Zhou17}.  Zhao studied $N$-derivations of Lie algebras in \cite{Zhao15}.  In this paper, we consider $n$-dervations of Lie color algebras and prove that every $n$-derivation of a perfect Lie color algebra with zero
center is a derivation.

\section{MAIN RESULTS}
\def\theequation{\arabic{section}. \arabic{equation}}
\setcounter{equation} {0}

Throughout the following section, $L$ always denotes a Lie color algebra over a
commutative ring $R$. A Lie color algebra $L$ is called perfect if the derived subalgebra
$[L, L]=L$. The center of $L$ is denoted by $Z(L)$. For a subset $S$ of $L$, denote by $C_L(S)$
the centralizer of $S$ in $L$. $Der(L)$ is the derivation algebra of $L$.  

\begin{definition}$^{\cite{Scheunert}}$
Let $\Gamma$ be an abelian group. A bicharacter on $\Gamma$ is a map $\epsilon: \Gamma\times \Gamma\rightarrow \mathbb{K}\backslash \{0\}$ satisfying
\begin{eqnarray*}
&&(1)~~ \epsilon(\a,\b)\epsilon(\b,\a)=1,\\
&&(2)~~\epsilon(\a,\b+\g)=\epsilon(\a,\b)\epsilon(\a,\g),\\
&&(3)~~ \epsilon(\a+\b,\g)=\epsilon(\a,\g)\epsilon(\b,\g),
\end{eqnarray*}
for any  $\a,\b\in \Gamma$.
\end{definition}

\begin{definition}$^{\cite{Scheunert}}$  A Lie color algebra  is a triple $(L, [\c,\c],  \epsilon)$ consisting of a $\Gamma$-graded space $L$,  a bilinear mapping $[\c,\c]: L\times L\rightarrow L$,  and a bicharacter $\epsilon$ on $\Gamma$   satisfying the following conditions,
\begin{eqnarray*}
&&[L_{x},L_{y}]\subseteq L_{xy},\\
&&[x,y]=-\epsilon(x, y)[y,x],\\
&&\epsilon(z, x)[x,[y,z]]+\epsilon(x, y)[y,[z,x]]+\epsilon(y, z)[z,[x,y]]=0,
\end{eqnarray*}
for any homogeneous elements $x, y, z\in L$.
\end{definition}

The $n$-derivation of $L$ is defined as follows:

\begin{definition}
An endomorphism of $R$-module $D$ of $L$ is called an $n$-derivation of
$L$.   For any $x_1, x_2, x_3,..., x_n\in L$, $D$  satisfies
\begin{eqnarray*}
&&D([[\c\c\c[[x_1,x_2],x_3],\c\c\c x_{n-1}],x_n])\\
&=& [[\c\c\c[[D(x_1),x_2],x_3],\c\c\c x_{n-1}],x_n]+\\
&& \sum_{i=2}^{n}\epsilon(D, \sum_{k=1}^{i-1}x_k)[[\c\c\c[\c\c\c[[x_1, x_2], x_3],\c\c\c D(x_i)],\c\c\c x_{n-1}], x_n].
\end{eqnarray*}
\end{definition}
Denote by $nDer(L)$ as the $R$ submodule spanned by all
$n$ derivations of $L$.
\smallskip

The main result of this article is the following theorem.

\begin{theorem}
Let $L$ be a Lie color algebra over $R$. If $\frac{1}{n-1}\in R$, $L$ is perfect and
has zero center, then we have that:

(1) ~$nDer(L)=Der(L).$

(2)~$nDer(Der(L))=ad(Der(L))$.
\end{theorem}

We proceed to prove the theorem in following lemmas.
 \begin{lemma}
 For any Lie color algebra $L$, $nDer(L)$ is closed under the usual Lie bracket. Furthermore, $nDer(L)$ is a Lie color algebra.
 \end{lemma}
 {\bf Proof.} For any $D_1,D_2\in nDer(L), x_1,x_2,...,x_n\in L$, we have
 \begin{eqnarray*}
&&D_1D_2([[\c\c\c[[x_1,x_2],x_3],\c\c\c x_{n-1}],x_n])\\
&=&D_1( [[\c\c\c[[D_2(x_1),x_2],x_3],\c\c\c x_{n-1}],x_n]+\\
&& \sum_{i=2}^{n}\epsilon(D_2, \sum_{k=1}^{i-1}x_k)[[\c\c\c[\c\c\c[[x_1, x_2], x_3],\c\c\c D_2(x_i)],\c\c\c x_{n-1}], x_n])\\
&=& [[\c\c\c[[D_1D_2(x_1),x_2],x_3],\c\c\c x_{n-1}],x_n]+\\
&& \sum_{i=2}^{n} \epsilon(D_1+D_2, \sum_{k=1}^{i-1}x_k)([[\c\c\c[\c\c\c[[x_1, x_2], x_3],\c\c\c D_1D_2(x_i)],\c\c\c x_{n-1}], x_n])\\
&& +\sum_{1\leq p< q\leq n}\epsilon(D_1, \sum_{k=1}^{p-1}x_k)\epsilon(D_2, \sum_{k=1}^{q-1}x_k)([\c\c\c[\c\c\c[\c\c\c D_1(x_p)], \c\c\c D_2(x_q)], \c\c\c x_n])\\
&& +\sum_{1\leq s< t\leq n} \epsilon(D_1,  \sum_{k=1}^{t-1}x_k+D_2)\epsilon(D_2, \sum_{k=1}^{s-1}x_k)([\c\c\c[\c\c\c[\c\c\c D_2(x_n)], \c\c\c D_1(x_t)], \c\c\c x_n]).
 \end{eqnarray*}
Similarly, we have
\begin{eqnarray*}
&&D_2D_1([[\c\c\c[[x_1,x_2],x_3],\c\c\c x_{n-1}],x_n])\\
&=&D_2( [[\c\c\c[[D_1(x_1),x_2],x_3],\c\c\c x_{n-1}],x_n]+\\
&& \sum_{i=2}^{n}\epsilon(D_1, \sum_{k=1}^{i-1}x_k)[[\c\c\c[\c\c\c[[x_1, x_2], x_3],\c\c\c D_1(x_i)],\c\c\c x_{n-1}], x_n])\\
&=& [[\c\c\c[[D_2D_1(x_1),x_2],x_3],\c\c\c x_{n-1}],x_n]+\\
&& \sum_{i=2}^{n} \epsilon(D_1+D_2, \sum_{k=1}^{i-1}x_k)([[\c\c\c[\c\c\c[[x_1, x_2], x_3],\c\c\c D_2D_1(x_i)],\c\c\c x_{n-1}], x_n])\\
&& +\sum_{1\leq p< q\leq n}\epsilon(D_1, \sum_{k=1}^{p-1}x_k)\epsilon(D_2, \sum_{k=1}^{q-1}x_k)([\c\c\c[\c\c\c[\c\c\c D_1(x_p)], \c\c\c D_2(x_q)], \c\c\c x_n])\\
&& +\sum_{1\leq s< t\leq n} \epsilon(D_1,  \sum_{k=1}^{t-1}x_k+D_2)\epsilon(D_2, \sum_{k=1}^{s-1}x_k)([\c\c\c[\c\c\c[\c\c\c D_2(x_n)], \c\c\c D_1(x_t)], \c\c\c x_n]).
 \end{eqnarray*}
Then, we have
\begin{eqnarray*}
&&[D_1,D_2]([[\c\c\c[[x_1,x_2],x_3],\c\c\c x_{n-1}],x_n])\\
&=& D_1D_2([[\c\c\c[[x_1,x_2],x_3],\c\c\c x_{n-1}],x_n])-\epsilon(D_1,D_2)D_2D_1([[\c\c\c[[x_1,x_2],x_3],\c\c\c x_{n-1}],x_n])\\
&=& [[\c\c\c[[[D_1,D_2](x_1),x_2],x_3],\c\c\c x_{n-1}],x_n] \\
&&+\sum_{i=2}^{n}\epsilon(D_1+D_2, \sum_{k=1}^{i-1}x_k)[[\c\c\c[\c\c\c[[x_1, x_2], x_3],\c\c\c [D_1, D_2](x_i)],\c\c\c x_{n-1}], x_n]).
\end{eqnarray*}
Hence, $[D_1,D_2]\in nDer(L)$. And the lemma is proved. \hfill $\square$
 \begin{lemma}
 If $L$ is a perfect Lie color algebra, then $ad(L)$ is an ideal of the Lie color algebra $nDer(L)$.
 \end{lemma}
{\bf Proof.} Let $D\in nDer(L), x\in L$. Since $L$ is perfect, there exists some finite index set
$I$ and $x_{ij}\in L$ such that
\begin{eqnarray*}
x=\sum_{i\in I}[\c\c\c[[x_{i1}, x_{i2}],x_{i3}],\c\c\c x_{in-1}].
\end{eqnarray*}
For any $z\in L$, we have
\begin{eqnarray*}
&&[D, adx](z)\\
&=& Dadx-\epsilon(D, \sum_{k=1}^{n-1}x_{ik})adx(D(z))\\
&=& D(x, z)-\epsilon(D, \sum_{k=1}^{n-1}x_{ik})[x, D(z)]\\
&=& D([\sum_{i\in I}[\c\c\c[[x_{i1}, x_{i2}],x_{i3}],\c\c\c x_{in-1}], z])-\epsilon(D, \sum_{k=1}^{n-1}x_{ik}) [\sum_{i\in I}[\c\c\c[[x_{i1}, x_{i2}],x_{i3}],\c\c\c x_{in-1}], D(z)]\\
&=&  \sum_{i\in I} D([[\c\c\c[[x_{i1}, x_{i2}],x_{i3}],\c\c\c x_{in-1}], z])- \sum_{i\in I}\epsilon(D, \sum_{k=1}^{n-1}x_{ik}) [[\c\c\c[[x_{i1}, x_{i2}],x_{i3}],\c\c\c x_{in-1}], D(z)]\\
&=&  \sum_{i\in I} [[\c\c\c[[D(x_{i1}), x_{i2}],x_{i3}],\c\c\c x_{in-1}], z]+ \sum_{i\in I} \epsilon(D, x_{i1})[[\c\c\c[[x_{i1}, D(x_{i2})],x_{i3}],\c\c\c x_{in-1}], z]\\
&& +\sum_{i\in I} \epsilon(D, x_{i1}+x_{i2}) [[\c\c\c[[x_{i1}, x_{i2}],D(x_{i3})],\c\c\c x_{in-1}], z]\\
&& +\c\c\c+\sum_{i\in I} \epsilon(D, x_{i1}+x_{i2}+x_{in-2}) [[\c\c\c[[x_{i1}, x_{i2}],x_{i3}],\c\c\c D(x_{in-1})], z]\\
&& +\sum_{i\in I}\epsilon(D, \sum_{k=1}^{n-1}x_{ik}) [[\c\c\c[[x_{i1}, x_{i2}],x_{i3}],\c\c\c x_{in-1}], D(z)] \\
&&-\sum_{i\in I}\epsilon(D, \sum_{k=1}^{n-1}x_{ik}) [[\c\c\c[[x_{i1}, x_{i2}],x_{i3}],\c\c\c x_{in-1}], D(z)]\\
&=& ad( \sum_{i\in I} [[\c\c\c[[D(x_{i1}), x_{i2}],x_{i3}],\c\c\c x_{in-1}], z]+ \sum_{i\in I} \epsilon(D, x_{i1})[[\c\c\c[[x_{i1}, D(x_{i2})],x_{i3}],\c\c\c x_{in-1}], z]+\\
&& \sum_{i\in I} \epsilon(D, x_{i1}+x_{i2}) [[\c\c\c[[x_{i1}, x_{i2}],D(x_{i3})],\c\c\c x_{in-1}], z]+\c\c\c+\\
&& \sum_{i\in I} \epsilon(D, x_{i1}+x_{i2}+x_{in-2}) [[\c\c\c[[x_{i1}, x_{i2}],x_{i3}],\c\c\c D(x_{in-1})], z])(z).
\end{eqnarray*}

By the arbitrariness of $z$, $[D, adx]$ is an inner derivation. Hence, $ad(L)$ is an ideal of $nDer(L)$. The proof is finished. \hfill $\square$

  \begin{lemma}
 If $L$ is a perfect Lie color algebra with zero center, then there exits an $R$-module
homomorphism $\delta: nDer(L)\rightarrow End(L), \delta(D)=\delta_D$  such that for all $x \in L$,
$D\in nDer(L)$, one has $[D, adx]=ad\delta_D(x)$.
 \end{lemma}
{\bf Proof.} By the proof of Lemma 2.6, if $L$ is perfect and has zero center, $D\in nDer(L)$, we can define a module endomorphism $\delta_D$ on $L$, such that
 \begin{eqnarray*}
x=\sum_{i\in I}[\c\c\c[[x_{i1}, x_{i2}],x_{i3}],\c\c\c x_{in-1}].
\end{eqnarray*}
Then we have
\begin{eqnarray*}
\delta_D(x)&=&\sum_{i\in I} [\c\c\c[[D(x_{i1}), x_{i2}],x_{i3}],\c\c\c x_{in-1}]+\\
&& \sum_{i\in I} \epsilon(D, x_{i1})[\c\c\c[[x_{i1}, D(x_{i2})],x_{i3}],\c\c\c x_{in-1}]+\\
&& \sum_{i\in I} \epsilon(D, x_{i1}+x_{i2}) [\c\c\c[[x_{i1}, x_{i2}],D(x_{i3})],\c\c\c x_{in-1}]+\c\c\c+\\
&& \sum_{i\in I} \epsilon(D, x_{i1}+x_{i2}+x_{in-2}) [\c\c\c[[x_{i1}, x_{i2}],x_{i3}],\c\c\c D(x_{in-1})].
\end{eqnarray*}
It is easy to check that  the definition is independent of the form of expression of $x$.  Hence, $\delta_D$ is well-defined and we have $[D, adx]=ad\delta_D(x)$, as desired.  \hfill $\square$
\medskip

By the map $\delta_D$ and the proof of Lemma 2.6, we have the following
lemmas.

 \begin{lemma}
 If $L$ is a perfect Lie color algebra with zero center, then  for all $D\in nDer(L)$, $\delta_D\in (n-1)Der(L)$.
 \end{lemma}
{\bf Proof.} For any $D\in nDer(L)$ and $x_1,x_2,\c\c\c, x_n\in L$, by Lemma 2.7 we have
\begin{eqnarray*}
[D, ad([\c\c\c[[x_1,x_2],x_3],\c\c\c x_{n-1}])]=ad\delta_D([\c\c\c[[x_1,x_2],x_3],\c\c\c x_{n-1}]).
\end{eqnarray*}
On the other hand, we have
\begin{eqnarray*}
&& [D, ad([\c\c\c[[x_1,x_2],x_3],\c\c\c x_{n-1}])]\\
&=& [D, [\c\c\c[[adx_1, adx_2], adx_3],\c\c\c adx_{n-1}]]\\
&=&[[D, [\c\c\c[[adx_1, adx_2], adx_3],\c\c\c adx_{n-2}]], adx_{n-1}]\\
&&+  \epsilon(D, \sum_{k=1}^{n-2}x_k)[[\c\c\c[[adx_1, adx_2], adx_3],\c\c\c adx_{n-2}], [D, adx_{n-1}]]\\
&=&[[D, [\c\c\c[[adx_1, adx_2], adx_3],\c\c\c adx_{n-2}]], adx_{n-1}]\\
&&+  \epsilon(D, \sum_{k=1}^{n-2}x_k)[[\c\c\c[[adx_1, adx_2], adx_3],\c\c\c adx_{n-2}], ad\delta_D(x_{n-1})]\\
&=& [[\c\c\c[[ad\delta_D(x_{1}), adx_2], adx_3],\c\c\c adx_{n-2}], adx_{n-1}]
\end{eqnarray*}
\begin{eqnarray*}
&&+ \epsilon(D,x_1) [[\c\c\c[[adx_{1}, ad\delta_D(x_2)],adx_3],\c\c\c adx_{n-2}], adx_{n-1}]\\
&&+\c\c\c+\epsilon(D, \sum_{k=1}^{n-2}x_k)[[\c\c\c[[adx_1, adx_2], adx_3],\c\c\c adx_{n-2}], [D, adx_{n-1}]]\\
&=& ad([[\c\c\c[[\delta_D(x_{1}), x_2], x_3],\c\c\c x_{n-2}], x_{n-1}]\\
&&+ \epsilon(D,x_1) [[\c\c\c[[x_{1}, \delta_D(x_2)],x_3],\c\c\c x_{n-2}], x_{n-1}]\\
&&+\c\c\c+\epsilon(D, \sum_{k=1}^{n-2}x_k)[[\c\c\c[[x_1, x_2], x_3],\c\c\c x_{n-2}], \delta_D(x_{n-1})]).
\end{eqnarray*}
Hence, we have
\begin{eqnarray*}
&&ad\delta_D([\c\c\c[[x_1,x_2],x_3],\c\c\c x_{n-1}])\\
&=& ad([[\c\c\c[[\delta_D(x_{1}), x_2], x_3],\c\c\c x_{n-2}], x_{n-1}]\\
&&+ \epsilon(D,x_1) [[\c\c\c[[x_{1}, \delta_D(x_2)],x_3],\c\c\c x_{n-2}], x_{n-1}]\\
&&+\c\c\c+\epsilon(D, \sum_{k=1}^{n-2}x_k)[[\c\c\c[[x_1, x_2], x_3],\c\c\c x_{n-2}], \delta_D(x_{n-1})]).
\end{eqnarray*}
Since $Z(L)=0$, then
\begin{eqnarray*}
&&\delta_D([\c\c\c[[x_1,x_2],x_3],\c\c\c x_{n-1}])\\
&=& [[\c\c\c[[\delta_D(x_{1}), x_2], x_3],\c\c\c x_{n-2}], x_{n-1}]\\
&&+ \epsilon(D,x_1) [[\c\c\c[[x_{1}, \delta_D(x_2)],x_3],\c\c\c x_{n-2}], x_{n-1}]\\
&&+\c\c\c+\epsilon(D, \sum_{k=1}^{n-2}x_k)[[\c\c\c[[x_1, x_2], x_3],\c\c\c x_{n-2}], \delta_D(x_{n-1})].
\end{eqnarray*}
By the arbitrariness of $x_1,x_2\c\c\c x_{n-2}, x_{n-1}$, we have  $\delta_D\in (n-1)Der(L)$, as required. \hfill $\square$
 \begin{lemma}
 If the base ring $R$ contains $\frac{1}{n-1}$, $L$ is perfect, then the centralizer of $ad(L)$ in
$nDer(L)$ is trivial, i.e., $C_{nDer(L)}(ad(L))=0$. In particular, the center of $nDer(L)$ is zero.
 \end{lemma}
{\bf Proof.}  Let $D\in C_{nDer(L)}(ad(L))$.  Then for any $x\in L$, $[D, adx]=0$.  Hence, for any $x,y\in L$, we have
\begin{eqnarray*}
D([x,y])-\epsilon(D, x)[x, D(y)]=[D, adx](y)=0.
\end{eqnarray*}
It follows that $
      D([x,y])=\epsilon(D, x)[x, D(y)].
      $
      Moreover, we have
      \begin{eqnarray*}
     D([x,y])= -\epsilon(x, y)D([y ,x])
     = -\epsilon(x, y) \epsilon(D, y) [y, D(x)]
     =[D(x), y].
      \end{eqnarray*}
      Hence, we have
     $
    D([x,y])=[D(x), y]=\epsilon(D, x)[x, D(y)].
      $

  For any $x_1,x_2,\c\c\c, x_{n-1}, x_n\in L$, we have
  \begin{eqnarray*}
&&D([[\c\c\c[[x_1,x_2],x_3],\c\c\c x_{n-1}],x_n])\\
&=& [[\c\c\c[[D(x_1),x_2],x_3],\c\c\c x_{n-1}],x_n]+\\
&& \sum_{i=2}^{n}\epsilon(D, \sum_{k=1}^{i-1}x_k)[[\c\c\c[\c\c\c[[x_1, x_2], x_3],\c\c\c D(x_i)],\c\c\c x_{n-1}], x_n]\\
&=& n D([[\c\c\c[[x_1,x_2],x_3],\c\c\c x_{n-1}],x_n]).
  \end{eqnarray*}
It follows that
        \begin{eqnarray*}
(n-1) D([[\c\c\c[[x_1,x_2],x_3],\c\c\c x_{n-1}],x_n])=0.
  \end{eqnarray*}
Since $\frac{1}{n-1}\in R$, we have

      \begin{eqnarray*}
D([[\c\c\c[[x_1,x_2],x_3]\c\c\c x_{n-1}],x_n])=0.
  \end{eqnarray*}
  Since $L$ is perfect, every element of $L$ can be expressed as the linear combination of
elements of the form $[x_1, [x_2, \c\c\c[x_{n-1}, x_n]\c\c\c]]$, we have that $D=0$.   \hfill $\square$
\medskip

The following lemma is easy to be proved.

\begin{lemma}
 For any Lie color algebra $L$, if $x\in D, D\in nDer(L)$, then $[D, adx]=ad(D(x))$.
\end{lemma}

Now we can prove the first conclusion of the theorem.

\begin{lemma}
If the base ring $R$ contains $\frac{1}{n-1}$, $L$ is perfect and has trivial center, then  $nDer(L)=Der(L)$.
\end{lemma}
{\bf Proof.} We apply a mathematical Induction.

(1) It is true for $n=3$ in \cite{Zhou13}.

(2) Suppose that it is true for $n=k$, we consider $n=k+1$.
Let $D$ be a $k+1$-derivation. By Lemma 2.7, we have $[D,adx]=ad\delta_D(x)$ for any $x\in L$,
where $\delta_D$  a $k$-derivation, then $\delta_D$ is a  derivation.
By Lemma 2.10, we have $[\delta_D, adx]=ad\delta_D(x)$.
 Moreover, we have $[D-\delta_D, adx]=0$ for any $x\in L$.
Hence, $D-\delta_D\in C_{nDer(L)}(ad(L))$. By Lemma 2.7, we have $D=\delta_D$, and we have that $k+1$-derivation is a  derivation.   Therefore,  $nDer(L)=Der(L)$.  \hfill $\square$
\begin{lemma}
If $L$ is a perfect Lie color algebra, $D\in nDer(nDer(L))$,  then  $D(ad(L))\subseteq ad(L)$.
\end{lemma}
{\bf Proof.} Since  $L$ is perfect, for any $x\in L$, there exist $x_{ij}\in L$  such that
 \begin{eqnarray*}
x=\sum_{i\in I}[[\c\c\c[[x_{i1}, x_{i2}],x_{i3}],\c\c\c x_{in-1}], x_{in}].
\end{eqnarray*}
Therefore, we have
\begin{eqnarray*}
&&D(adx)\\
&=& \sum_{i\in I}D(ad[[\c\c\c[[x_{i1}, x_{i2}],x_{i3}],\c\c\c x_{in-1}], x_{in}])\\
&=& \sum_{i\in I}D([[\c\c\c[[adx_{i1}, adx_{i2}],adx_{i3}],\c\c\c adx_{in-1}], adx_{in}])\\
&=& \sum_{i\in I} ([[\c\c\c[[D(adx_{i1}), adx_{i2}],adx_{i3}],\c\c\c adx_{in-1}], adx_{in}]\\
&&+\epsilon(D, \sum_{k=1}^{l-1}x_{ik}[[\c\c\c[[adx_{i1}), adx_{i2}],adx_{i3}],\c\c\c D( adx_{il-1})],\c\c\c adx_{in-1}], adx_{in}]).
\end{eqnarray*}
Since $D\in nDer(nDer(L))$, then $D(adx_{im})\subseteq nDer(L)$. By Lemma 2.6, we have  $D(ad(L))\subseteq ad(L)$. \hfill $\square$
\begin{lemma}
Suppose that $L$ is a perfect Lie color algebra with zero center, $D\in nDer(nDer(L))$. If  $D(ad(L))=0$, then $D=0$.
\end{lemma}
{\bf Proof.} Since  $L$ is perfect, for any $x\in L$, there exist $x_{ij}\in L$  such that
 \begin{eqnarray*}
x=\sum_{i\in I}[\c\c\c[[x_{i1}, x_{i2}],x_{i3}],\c\c\c x_{in-1}].
\end{eqnarray*}
For any $d\in nDer(L)$, we have
\begin{eqnarray*}
&&[adx, D(d)]\\
&=& [\sum_{i\in I}[\c\c\c[[ad_{i1},ad_{i2}], ad_{i3}],\c\c\c ad_{in-1}],  D(d)]\\
&=& \sum_{i\in I}(\epsilon(D, \sum_{k=1}^{n-1}x_{ik})D([[\c\c\c[[adx_{i1}, adx_{i2}],adx_{i3}],\c\c\c adx_{in-1}], d])\\
&& -\epsilon(D, \sum_{k=1}^{n-1}x_{ik}) [[\c\c\c[[D(adx_{i1}), adx_{i2}],adx_{i3}],\c\c\c adx_{in-1}], d]\\
&&-\sum_{l=2}^{n-1}\epsilon(D, \sum_{k=1}^{n-1}x_{ik}+\sum_{s=1}^{l-1}x_{ip})[[\c\c\c[\c\c\c[[adx_{i1}, adx_{i2}],adx_{i3}],\c\c\c D(adx_{il})], \c\c\c adx_{in-1}], adx_{in}], d]).
\end{eqnarray*}
By Lemma 2.6, we have
\begin{eqnarray*}
[[\c\c\c[[adx_{i1}, adx_{i2}],adx_{i3}],\c\c\c adx_{in-1}], d]\in ad(L).
\end{eqnarray*}
Then
\begin{eqnarray*}
 D([[\c\c\c[[adx_{i1}, adx_{i2}],adx_{i3}],\c\c\c adx_{in-1}], d])=0.
\end{eqnarray*}
It follows that $
       D([adx, d])=0.
      $
       Moreover, we have
       $
       D(d)\in C_{nDer(L)}(ad(L)).
      $
       By Lemma 2.9, we have $D(d)=0$ and therefore $D=0$, as desired.  \hfill $\square$

\begin{lemma}Let $L$ be a Lie color algebra over $R$. Suppose that $\frac{1}{n-1}\in R$, $L$
is perfect and has zero center. If $D\in nDer(nDer(L))$, then there exists $d\in Der(L)$ such
that for any $x\in L$, $D(adx)=ad(d(x))$.
\end{lemma}
{\bf Proof.} For any $D\in nDer(nDer(L))$ and $x\in L$. By Lemma 2.12, $D(adx)\in ad(L)$. Let $y\in L$
and $D(adx)=ady$. Since the center $Z(L)$ is trivial, such $y$ is unique. Clearly, the
map  $d: x\rightarrow y$  is a $R$-module endomorphism of $L$.

For any $x_1,x_2,\c\c\c, x_n\in L$, we have
\begin{eqnarray*}
&& ad(d([[\c\c\c[[x_1,x_2],x_3],\c\c\c x_{n-1}],x_n]))\\
&=& D(ad([[\c\c\c[[x_1,x_2],x_3],\c\c\c x_{n-1}],x_n]))\\
&=& D([[\c\c\c[[adx_1,adx_2],adx_3],\c\c\c adx_{n-1}],adx_n])\\
&=& [[\c\c\c[[D(adx_1), adx_2],adx_3],\c\c\c adx_{n-1}],adx_n]\\
&& \epsilon(D, \sum_{k=1}^{i-1}x_k[[\c\c\c[\c\c\c[[adx_1, adx_2],adx_3]\c\c\c D(adx_i)],\c\c\c adx_{n-1}],adx_n]\\
&=& [[\c\c\c[[ad(d(x_1)), adx_2],adx_3],\c\c\c adx_{n-1}],adx_n]\\
&& \epsilon(D, \sum_{k=1}^{i-1}x_k[[\c\c\c[\c\c\c[[adx_1, adx_2],adx_3],\c\c\c ad(d(x_i))],\c\c\c adx_{n-1}],adx_n]\\
&=& ad([[\c\c\c[[d(x_1), x_2],x_3],\c\c\c x_{n-1}],x_n]\\
&& \epsilon(D, \sum_{k=1}^{i-1}x_k[[\c\c\c[\c\c\c[[x_1, x_2],x_3],\c\c\c d(x_i)],\c\c\c x_{n-1}],x_n]).
\end{eqnarray*}
Since $Z(L)=0$, we have
\begin{eqnarray*}
&&d([[\c\c\c[[x_1,x_2],x_3],\c\c\c x_{n-1}],x_n])\\
&=&[[\c\c\c[[d(x_1), x_2],x_3],\c\c\c x_{n-1}],x_n]\\
&& \epsilon(D, \sum_{k=1}^{i-1}x_k[[\c\c\c[\c\c\c[[x_1, x_2],x_3],\c\c\c d(x_i)],\c\c\c x_{n-1}],x_n].
\end{eqnarray*}
Thus, $d\in nDer(L)$.  By Lemma 2.11, we have  $d\in Der(L)$ and the lemma is finished. \hfill $\square$
\begin{lemma}
Let $L$ be a Lie color algebra over $R$. If $\frac{1}{n-1}\in R$, $L$ is perfect and
has zero center, then $nDer(Der(L))=ad(Der(L))$.
\end{lemma}
{\bf Proof.} For any $D\in nDer(nDer(L))$ and $x\in L$,  there exists $d\in Der(L)$ such
that for any $x\in L$, $D(adx)=ad(d(x))$. By Lemma 2.10, we have $ad(d(x))=[d, adx]$. Hence, we have
\begin{eqnarray*}
 D(adx)=ad(d(x))=[d, adx]=ad(d)(adx).
\end{eqnarray*}
Thus $
 D-ad(d)(adx)=0.
$
By Lemma 2.13, $D=ad(d)$. Therefore, $nDer(Der(L))=ad(Der(L))$.   \hfill $\square$

\begin{corollary}
Let $L$ be a Lie superalgebra over  $R$. If $\frac{1}{n-1}\in R$, $L$ is perfect and
has zero center, then we have that:

(1)~ $nDer(L)=Der(L)$

(2)~$nDer(Der(L))=ad(Der(L))$.
\end{corollary}

\begin{corollary}
Let $L$ be a Lie algebra over  $R$. If $\frac{1}{n-1}\in R$, $L$ is perfect and
has zero center, then we have that:

(1)~ $nDer(L)=Der(L)$

(2) ~$nDer(Der(L))=ad(Der(L))$.
\end{corollary}

 \begin{center}
 {\bf ACKNOWLEDGEMENT}
 \end{center}

  The paper is supported by the Youth Project for Natural Science Foundation of Guizhou provincial department of education (No. KY[2018]155),
 the innovative exploration and academic seedling project of Guizhou University of Finance and Economics (No. [2017]5736-023).


\begin{thebibliography}{99}



\bibitem{Chen13} L. Chen, Y. Ma, L. Ni.    Generalized derivations of Lie color algebras.  Results in Mathematics,   63 (3-4)(2013), 923-936.

\bibitem{Lian16} H. Lian, C. Chen.   $N$-derivations for finitely generated graded Lie algebras.  Algebra Colloquium,  23(2016), 205-212.

\bibitem{Scheunert} M. Scheunert. Generalized Lie algebras. Journal of Mathematical Physics, 20(1979), 712-720.

\bibitem{Sun18} B. Sun, L.  Chen.   Double Derivations of $n$-Lie Superalgebras.  Algebra Colloquium,  25(2018), 161-180.

    \bibitem{Wang13}  Y. Wang, Y. Wang, Y. Du.  $n$-Derivations of triangular algebras.   Linear Algebra and Its Applications, 439 (2) (2013), 463-471.


  \bibitem{Zhao15} D. Zhao.  $N$-derivations of Lie algebras.    Master thesis of Southeast University  (China). 2015, 57 pp.

\bibitem{Zhou17}  J. Zhou, L. Chen, Y. Ma.  Triple derivations and triple homomorphisms of perfect Lie superalgebras.  Indagationes Mathematicae, 28 (2017),  436-445.

\bibitem{Zhou13}   J. H. Zhou.    Triple derivations of perfect Lie algebras. Communications in Algebra,    41 (2013), 1647-1654.

\end{thebibliography}
\end{document}